\documentclass[fleqn]{mat01}
\usepackage{times,mathtimy,amssymb,latexsym}%,krepsf,rangecite}
\begin{document}

\setcounter{page}{293} \firstpage{293}

\font\xxxx=msam10 at 10pt
\def\epr{\mbox{\xxxx{\char'245\ \,}}}

\newtheorem{coro}[defin]{\rm COROLLARY}
\newtheorem{exam}[defin]{\it Example}

\title{Pro-torus actions on Poincar\'{e} duality spaces}

\markboth{Ali \"{O}zkurt and Do\u{g}an D\"{o}nmez}{Pro-torus
actions on Poincar\'{e} duality spaces}

\author{ALI \"{O}ZKURT and DO\u{G}AN D\"{O}NMEZ}

\address{Department of Mathematics, \c{C}ukurova University, 01330-Adana,
Turkey\\
\noindent E-mail: aozkurt@mail.cu.edu.tr; donmez@mail.cu.edu.tr}

\volume{116}

\mon{August}

\parts{3}

\pubyear{2006}

\Date{MS received 3 October 2005; revised 11 May 2006}

\begin{abstract}
In this paper, it is shown that some of the results of torus
actions on Poincar\'{e} duality spaces, Borel's dimension formula
and topological splitting principle to local weights, hold if
`torus' is replaced by `pro-torus'.
\end{abstract}


\maketitle

\section{Introduction}

In the theory of linear representations of compact connected Lie
groups, the crucial first step is restriction to the maximal
torus. Hsiang \cite{10} has suggested that the study of
topological transformation groups proceed in the same way, and has
offered the concepts of local weights and $F$-varieties to
generalize the linear notions of weights and weight\break spaces.

It is well-known that locally compact groups can be approximated
by Lie groups. This means if $G$ is a locally compact group, with
finitely many components, then $G$ has arbitrarily small compact
normal subgroup $N$ such that $G/N$ is a Lie group.

$G$ is called a $k$ dimensional pro-torus (i.e compact, connected,
abelian group) if $G$ has a totally disconnected closed subgroup
$N$ such that $G/N \simeq T^{k}$, a $k$ torus. Furthermore, the
dimension of $G$ is infinite if $G$ has a totally disconnected
closed subgroup $N$ such that $G/N$ is a infinite dimensional
torus. If a $k$ dimensional pro-torus $G$ acts effectively on a
Hausdorff space $X$, (all actions are assumed to be continuous),
then there is an induced action of the $k$-torus $G/N$ on the
orbit space $X/N$ and $X\rightarrow\ X/N$ induces a homeomorphism
$X^{G}\approx(\ X/N)^{\ G/N}$. (Here $X^{G}$ denotes fixed point
set of action $G$.) The orbit space of the action of $G/N$ on $X/N
$ is homeomorphic to the orbit space $X/G$. On the other hand, the
orbit space $X/N$ inherits global cohomological properties from
the space $X$. Thus many questions about the cohomological
properties of orbit spaces and fixed point sets of actions of
pro-torus are reduced to questions about torus.

Let $X$ be a connected Hausdorff topological space such that
$H^{*}(X,\mathbb{Q})$ finite dimensional. We say that $X$ is a
Poincar\'{e} duality space over the rationals of formal dimension
$n>0$, and we write $fd(X)=n$, if $H^{*}(X,\mathbb{Q})$ is a
finite dimensional vector space over $\mathbb{Q}$ and
$H^{*}(X,\mathbb{Q}) \otimes H^{n-*}(X,\mathbb{Q})\rightarrow
H^{n}(X,\mathbb{Q})=\mathbb{Q}$ is a nonsingular pairing. It is a
well-known result proven by Chang and Skjelberg \cite{7} and
Bredon \cite{4}, that each component of the fixed point set of a
torus action also satisfies Poincar\'{e} duality. Moreover, Borel
formula also holds. In this paper these results are generalized
for pro-torus actions using some techniques of Biller \cite{2}.
Biller showed that the fixed point set of pro-torus in a rational
cohomology manifold is a rational cohomology manifold for even
codimension.

From now on, $X$ will be a compact Poincar\'{e} duality space over
the rationals of formal dimension $n$, and all cohomologies are
sheaf cohomology over $\mathbb{Q}$.

\begin{theorem}[\cite{6,11}]
Let $N$ be a totally disconnected compact group which acts on a
locally compact space $X$. Then the orbit projection
$X\rightarrow\ X/N$ induces an isomorphism
\begin{equation*}
H_{c}^{*}(X/N,\mathbb{Q})\simeq H_{c}^{*}(X,\mathbb{Q})^{N}.
\end{equation*}
\vspace{-2pc}
\end{theorem}

\section{Components of fixed point sets and local weights}

In this section, we define some cohomology classes for each
components of fixed point set which we call local weights along each
component.

\setcounter{defin}{0}
\begin{theorem}[\!]
Let $G$ be a compact connected abelian group which acts on a
compact Poincar\'{e} duality space $X$ of formal dimension $n$.
Then each connected component of $X^{G}$ is also a Poincar\'{e}
duality space.
\end{theorem}

\begin{proof}
Let $G$ be a compact connected abelian group of finite dimension
which acts on a compact Poincar\'{e} duality space $X$ of formal
dimension $n$. Let $N$ be a  totally disconnected closed subgroup
of $G$ such that $G/N \simeq T^{k}$. ($k$ is called the dimension
of\break $G$.)

Since $G$ is connected, its action (hence that of $N$) on
$H^{*}(X)$ is trivial (see II.10.6, II.11.11 of \cite{5}). The
orbit space of $X/N$ is a Poincar\'{e} duality space such that
$fd(X/N)=fd(X)=n$ by Theorem~1.1. Let $F'\subset X^{G}$ be a
nonempty connected component of fixed point set $X^{G}$. Since the
space $(X/N)^{\ G/N}$ is homeomorphic to $X^{G}$, $F'$ can be
considered to be a connected component of $(X/N)^{\ G/N}$ by
identifying $F'$ with its homeomorphic image. So $F'$ is a
Poincar\'{e} duality space such that $fd(F')=r$.

If the dimension of $G$ is infinite then we find a totally
disconnected closed subgroup $N$ of $G$ such that $G/N$ is an
infinite dimensional torus and every sub-torus of $G/N$ is of the
form $HN/N$ for a pro-torus $H\leq G$ \cite{9}. Let $F'$ be a
connected component of $X^{G}$. Choose a finite dimensional
pro-torus $H_{1}\leq G$ whose action on $X$ is not trivial. Let
$F_{1}$ be the connected component of $X^{H_{1}}$ which contains
$F'$. If $F'\neq F_{1}$ then the induced action of $G/H_{1}$ on
$F_{1}$ is not trivial, hence we may choose a finite dimensional
pro-torus $H_{2}\leq G$ which contains $H_{1}$ and acts
non-trivially on $F_{1}$. So replacing $F_{1}$ and $H_{1}$ by $X$
and $G$ we can find a pro-torus $H_{2}$ such that $H_{1}\leq
H_{2}$ and $H_{2}$ acts non-trivially on $F_{1}$. Thus we have
constructed a properly descending sequence of Poincar\'{e} duality
spaces $F_{j}$ by induction. Since formal dimensions decrease
strictly, we reach a finite dimensional pro-torus $H_{m}\leq G$
such that $F'$ is a component of the fixed point set $X^{H_{m}}$.
Thus $F'$ is also a Poincar\'{e} duality
space.\hfill\epr\vspace{-.5pc}
\end{proof}

\begin{theorem}[\!]
Let $S$ be the set of all closed{\rm ,} connected{\rm ,}
codimension $1$ subgroups of $G$. Then the formula
\begin{equation*}
n-r = \sum_{H\in S} (n(H)-r)
\end{equation*}
holds where $n(H)$ denotes the formal dimension of connected
component of $X^{H}$ containing $F'$. We will denote this set by
$F'(H)$.
\end{theorem}

\begin{proof}
If the dimension of $G$ is infinite then $G$ has closed connected
subgroups of arbitrarily high finite dimension \cite{9}. We may
therefore assume that $G$ is finite dimensional. Choose $H\in S$.
Then fixed point set, $X^{H}$, is invariant under $G$. So
$H^{*}(X^{H})\simeq H^{*}(X^{H}/N)$ by Theorem $1.1$. Moreover
\hbox{$X^{H}\!/\!N=(X/N)^{H}=(X/N)^{HN/N}$}. Since
$H^{*}(X^{H})\simeq H^{*}((X/N)^{HN/N})$ there is a one-to-one
correspondence between connected components of $X^{H}$ and
$(X/N)^{HN/N}$. Let $F'\subset X^{G}$ be a connected component.
Then there exists $F'\subset F'(HN/N)$ such that $F'(HN/N)$ is a
connected component of $(X/N)^{HN/N}$. Moreover, every codimension
$1$ sub-torus of $G/N$ is of the form $HN/N$ for a unique $H\in S$
\cite{9}. So $F'(HN/N)$ is a Poincar\'{e} duality space such that
$fd(F'(HN/N))=n(HN/N)$.

Let $F'(H)$ be the connected component of $X^{H}$ containing $F'$
and assume it corresponds to $F'(HN/N)$ under the isomorphism
$H^{*}(X^{H})\simeq H^{*}((X/N)^{HN/N})$. Then it is clear that
$F'(H)$ is a Poincar\'{e} duality space and
$fd(F'(H))=fd(F'(HN/N))$. Finally, consider the action of $G/N
\simeq T^{k}$ on $X/N$ which is also a Poincar\'{e} duality space.
Then we can deduce the theorem from Borel's formula for torus
actions on Poincar\'{e} duality spaces:
\begin{align*}
n-r = \sum_{HN/N} (fd(F'(HN/N))-r) = \sum_{H\in S} (n(H)-r).
\end{align*}
$\left.\right.$\vspace{-1.5pc}

\hfill\epr
\end{proof}

If $n(H)>r$ then we will say that $H$ is a local weight at $F'$.

\begin{theorem}[\!]
$H$ is a local weight at $F'$ if and only if $i^{*}\hbox{\rm :}\
H^{*}(F'(H))\rightarrow H^{*}(F')$ {\rm (}induced by inclusion
$i\hbox{\rm :}\ F'\hookrightarrow F'(H)${\rm )} is not an
isomorphism.
\end{theorem}

\begin{proof}
If $H$ is a local weight at $F'$ then it is trivial that
$i^{*}\hbox{\rm :}\ H^{*}(F'(H))\rightarrow H^{*}(F')$ is not an
isomorphism. Conversely, assume that $H$ is not a local weight at
$F'$. Since $fd(F'(HN/N))=fd(F'(H))$, $HN/N$ is not a local weight
at $F'$ (considering the action of $G/N$ on $X/N$). Thus
$j^{*}\hbox{\rm :}\ H^{*}(F'(HN/N))\rightarrow H^{*}(F')$, where
$j\hbox{\rm :}\ F'\hookrightarrow F'(HN/N)$ is an isomorphism.
Since $H^{*}(F'(H))\simeq H^{*}(F'(HN/N))$, $i^{*}\hbox{\rm :}\
H^{*}(F'(H))\rightarrow H^{*}(F')$ is an isomorphism.\hfill\epr
\end{proof}

The Borel theory is reflected in the algebraic properties of the
equivariant cohomology ring $H_{G}^{*}(X)$ which is defined below.

Let $B_{G}$ be a classifying space of $G$ and $X_{G}=(X\times
E_{G})/G$ be the balanced product, where $E_{G}\rightarrow B_{G}$
is a universal principal $G$-bundle. Then the equivariant
cohomology ring of $X$ is $H_{G}^{*}(X)=H^{*}(X_{G})$. Let
$B_{\pi}\hbox{\rm :}\ B_{G}\rightarrow B_{G/N}$ be the mapping
induced by the canonical epimorphism $\pi\hbox{\rm :}\
G\rightarrow G/N$, since $G/N\simeq T^{k}$, $H^{*}(B_{G/N})$ is a
polynomial algebra over $\mathbb{Q}$. Let $R=\hbox{Im}\,
B_{\pi}^{*}$ where $B_{\pi}^{*}\hbox{\rm :}\ H^{*}(B_{G/\
N})\rightarrow H^{*}(B_{G})$. $R$ is a subring of $H^{*}(B_{G})$.

Let $f_{1}$ be a generator of $H^{r}(F')$ where $f_{1}\in
H^{*}(F')\subset H^{*}(X^{G})\subset
H_{G}^{*}(X^{G})=H^{*}(B_{G})\otimes H^{*}(X^{G})$. For $H \in S$,
let us define the ideals:
\begin{equation*}
I_{f_{1}}^{G}=\{a\in R\hbox{\rm :}\ a\otimes f_{1}\in
\hbox{Im}\,\{H_{G}^{*}(X)\rightarrow H_{G}^{*}(X^{G})\}\}.
\end{equation*}
and
\begin{equation*}
I_{f_{1}}^{H}=\{a\in R\hbox{\rm :}\ a\otimes f_{1}\in
\hbox{Im}\,\{H_{G}^{*}(X^{H})\rightarrow H_{G}^{*}(X^{G})\}\}.
\end{equation*}

We say that $X/N$ is totally nonhomologous to zero in
$(X/N)_{G/N}\rightarrow B_{G/N}$ with respect to rational
cohomology if
\begin{equation*}
j^{*}\hbox{\rm :}\ H_{G/N}^{*}(X/N)\rightarrow H^{*}(X/N)
\end{equation*}
is surjective. It is a well-known fact that
$\dim_{\mathbb{Q}}H^{*}((X/N)^{G/N})=\dim_{\mathbb{Q}}H^{*}(X/N)$
if and only if Leray--Serre spectral sequence associated to
$(X/N)_{G/N}\rightarrow B_{G/N}$ degenerates, i.e. $X/N$ is
totally non-homologous to zero in $(X/N)_{G/N}\rightarrow B_{G/N}$
with respect to rational cohomology.

\begin{theorem}[\!]
If $X/N$ is totally nonhomologous to zero in
$(X/N)_{G/N}\rightarrow B_{G/N}${\rm ,} then $I_{f_{1}}^{G}$ is a
principal ideal in $R${\rm ,} generated by an element $\xi\in
H^{n-r}(B_{G})$. Furthermore{\rm ,} $\xi$ splits as a product of
homogeneous linear factors in $H^{*}(B_{G})${\rm ,} such that each
linear factor corresponds to a local weight at $F'$ and has
multiplicity $\frac{n(H)-r}{2}$.
\end{theorem}

\begin{proof}
Let $I_{f_{1}}^{G/N}= \{a\in H^{*}(B_{G/N})\hbox{\rm :}\ a\otimes
f_{1}\in \hbox{Im}\,{H_{G/N}^{*}(X/N)\rightarrow
H_{G/N}^{*}(X/N)^{G/N}}\}$. If $a\in I_{f_{1}}^{G/N}$, then
$a\otimes f_{1}\in \hbox{Im}\,{H_{G/N}^{*}(X/N)\rightarrow
H_{G/N}^{*}(X/N)^{G/N}}$. Consider the commutative diagram
\begin{equation*}
\begin{array}{ccccccc}
&   & X & \rightarrow  & X/N \\
&  & \downarrow  &  & \downarrow &  &
\\
X^{G}\times B_{G} & \hookrightarrow  & X_{G} & \rightarrow  &
(X/N)_{G/N} & \hookleftarrow &(X/N)^{G/N}\times B_{G/N} \\
&  & \downarrow  &  & \downarrow &  & \\
&    & B_{G} & \rightarrow  & B_{G/N}\end{array}.
\end{equation*}
This diagram induces the following commutative diagram:
\begin{equation*}
\begin{array}{ccccccc}
&   & H^{*}(X/N)& \rightarrow  & H^{*}(X) \\
&  & \uparrow  &  & \uparrow &  & \\
H_{G/N}^{*}((X/N)^{G/N})
 & \leftarrow & H_{G/N}^{*}(X/N) &
\rightarrow & H_{G}^{*}(X) &
\rightarrow &H_{G}^{*}(X^{G}) \\
&  & \uparrow  &  & \uparrow &  & \\
&    & H^{*}(B_{G/N}) & \rightarrow  & H^{*}(B_{G})\end{array}.
\end{equation*}

It is clear that $B_{\pi}^{*}(a)\in
\hbox{Im}\,\{H_{G}^{*}(X)\rightarrow H_{G}^{*}(X^{G})\}$. On the
other hand, let $b=B_{\pi}^{*}(a)\in I_{f_{1}}^{G}$. Since $X/N$
is totally nonhomologous to zero in $(X/N)_{G/N}\rightarrow
B_{G/N}$, it is clear that $a\in I_{f_{1}}^{G/N}$. So restriction
of $B_{\pi}^{*}$ to $I_{f_{1}}^{G/N}$ induces an epimorphism
$I_{f_{1}}^{G/N}\rightarrow I_{f_{1}}^{G}$. Let $H$ be a local
weight at $F'\subset X^{G}$. Then $HN/N$ is a local weight at
$F'\subset (X/N)^{G/N}$. Let $w_{HN/N}\in H^{2}(B_{G/N})$ be the
corresponding cohomology class. Then
$w_{H}=B_{\pi}^{*}(w_{HN/N})\in H^{2}(B_{G})$ is the corresponding
cohomology class of $H$. It is well-known that $I_{f_{1}}^{G/N}$
is a principal ideal in $H^{*}(B_{G/N})$ generated by
${\prod}(w_{HN/N})^{m_{H}}\in H^{n-r}(B_{G/N})$ \cite{1}, where
$\frac{n(H)-r}{2}$. Then it is easy to see that $I_{f_{1}}^{G}$ is
a principal ideal in $R$ generated by ${\prod }(w_{H})^{m_{H}}$.
This finishes the proof.\hfill\epr
\end{proof}

\begin{coro}$\left.\right.$\vspace{.5pc}

\noindent If $H$ is a local weight at $F'$ then the ideal
$I_{f_{1}}^{H}$ is a principal maximal ideal with respect to the
property $I_{f_{1}}^{H}\neq R$ which is generated by $w_{H}$.
\end{coro}

\begin{proof}
If $H$ is a local weight at $F'$, then $HN/N$ is a local weight at
$F'$ (considering $G/N$ action on $X/N$). Then the ideal
\begin{align*}
I_{f_{1}}^{HN/N} &= \{a\in H^{*}(B_{G/N})\hbox{\rm :}\ a\otimes
f_{1}\in \hbox{Im}\,\{H_{G/N}^{*}((X/N)^{HN/N})\\[.4pc]
&\quad\,\rightarrow H_{G/N}^{*}((X/N)^{G/N})\}\}
\end{align*}
is a maximal ideal with respect to the property
$I_{f_{1}}^{HN/N}\neq H^{*}(B_{G/N})$ which is generated by
$w_{HN/N}$ \cite{1}. Let us replace $X^{H}$ by $X$ and consider
the $G/N$ action on $X^{H}/N=(X/N)^{HN/N}$. We have
$\dim_{\mathbb{Q}}H^{*}((X/N)^{G/N})=\dim_{\mathbb{Q}}H^{*}(X/N)$,
since $X/N$ is totally nonhomologous to zero in
$(X/N)_{G/N}\rightarrow B_{G/N}$. Let $G/N=S^{1}\times HN/N$. Then
\begin{align*}
\dim_{\mathbb{Q}}H^{*}((X/N)^{G/N}) &=
\dim_{\mathbb{Q}}H^{*}((X^{H}/N)^{S^{1}})\leq
\dim_{\mathbb{Q}}H^{*}(X^{H}/N)\\[.4pc]
&\leq \dim_{\mathbb{Q}}H^{*}(X/N).
\end{align*}
Thus
$\dim_{\mathbb{Q}}H^{*}(X^{H}/N)=\dim_{\mathbb{Q}}H^{*}(X/N)$.
This implies that $X^{H}/N$ is totally nonhomologous to zero in
$(X^{H}/N)_{G/N}\rightarrow B_{G/N}$. Consider the commutative
diagram:
\begin{equation*}
\begin{array}{ccccccc}
&   & X ^{H}& \rightarrow  & X^{H}/N \\
&  & \downarrow  &  & \downarrow &  & \\
X^{G}\times B_{G} & \hookrightarrow  & X^{H}_{G} & \rightarrow  &
(X^{H}/N)_{G/N} & \hookleftarrow &(X/N)^{G/N}\times B_{G/N} \\
&  & \downarrow  &  & \downarrow &  & \\
&    & B_{G} & \rightarrow  & B_{G/N}\end{array}.
\end{equation*}
It is clear that $B_{\pi}^{*}$ induces the epimorphism
$I_{f_{1}}^{HN/N}\rightarrow I_{f_{1}}^{H}$. This ends the
proof.\hfill\epr
\end{proof}

\section{Applications}

\setcounter{defin}{0}
\begin{exam}
{\rm Let $X$ be a $G$-space, $G$ finite-dimensional compact
connected abelian group and $H^{*}(X)=H^{*}(S^{n})$. Let $N$ be a
totally disconnected closed subgroup of $G$ such that $G/N\simeq
T^{k}$ is a torus group. Since the action of $G$, and hence that
of $N$, on $H^{*}(X)$ is trivial, $H^{*}(X/N)=H^{*}(S^{n})$ by
Theorem~$1.1$. So the space $X^{G}$, which is homeomorphic to
$(X/N)^{G/N}$, has the rational cohomology of $S^{r}$ for some
$r\in {-1,\dots,n}$ such that $n-r$ is even. Assume that
$X^{G}\neq \emptyset$. Therefore,
$\dim_{\mathbb{Q}}H^{*}(X^{G})=\dim_{\mathbb{Q}}H^{*}(X)=2$. Let
$T'$ be a geometric weight for the $G/N$ space $X/N$ (see
\cite{10}). Then $T'=HN/N$ for a unique $H\in S$. Since
$\dim_{\mathbb{Q}}H^{*}((X/N)^{G/N})=\dim_{\mathbb{Q}}H^{*}((X^{H})/N)$,
$X^{H}/N$ is totally nonhomologous to zero in
$(X^{H}/N)_{G/N}\rightarrow B_{G/N}$.}
\end{exam}

\begin{coro}$\left.\right.$\vspace{.5pc}

\noindent  Let $A$ be a closed invariant subspace of $X$ such that
inclusion $i\hbox{\rm :}\ A\hookrightarrow X$ induces an
isomorphism $H^{*}(X)\rightarrow H^{*}(A)$. Also assume that $X/N$
is totally nonhomologous to zero in $(X/N)_{G/N}\rightarrow
B_{G/N}$. Then the local weights of $A$ are equal to the local
weights of $X$.
\end{coro}

\begin{proof}
Recall that a space $X$ is called finitistic if every open
covering of $X$ has a finite dimensional open refinement. (The
dimension of a covering is the dimension of its nerve, which is
one less than the maximum number of members of the covering which
intersect nontrivially.) Clearly every compact space is
finitistic. Consider $G/N\simeq T^{k}$ space $X/N$ and its closed
invariant subspace $A/N$. It is easy to see that inclusion
$A/N\hookrightarrow X/N$ induces an isomorphism
$H^{*}(X/N)\rightarrow H^{*}(A/N)$ by Theorem~1.1. So $(X/N,A/N)$
is totally nonhomologous to zero in $(X/N,A/N)\rightarrow
((X/N)_{G/N},(A/N)_{G/N})\rightarrow B_{G/N}$. Thus inclusion
$(A/N)^{G/N}\hookrightarrow (X/N)^{G/N}$ induces isomorphism
$H^{*}((X/N)^{G/N})\simeq H^{*}((A/N)^{G/N})$. This is essentially
Theorem~1.6, ch.~VII of \cite{3}. The proof in Bredon is for the
case where $G = S^{1}$ or $Z_{p}$. But one gets the result for
higher-rank tori and $p$-tori using induction. We must assume
finitistic orbit space for $S^{1}$ action but this is now known by
Deo and Tripathi \cite{8}. Let us consider the Leray--Serre
spectral sequences of $A/N\rightarrow (A/N)_{G/N}\rightarrow
B_{G/N}$ and $X/N\rightarrow (X/N)_{G/N}\rightarrow B_{G/N}$. It
is clear that $H_{G/N}^{*}(X/N)\simeq H_{G/N}^{*}(A/N)$ by
Zeeman's comparison theorem. So $I_{f_{1}}^{G/N}= J_{f_{1}}^{G/N}$
where
\begin{align*}
J_{f_{1}}^{G/N}= \{a\in H^{*}(B_{G/N})\hbox{\rm :}\ a\otimes
f_{1}\in \hbox{Im}\,{H_{G/N}^{*}(A/N)\rightarrow
H_{G/N}^{*}(A/N)^{G/N}}\}.
\end{align*}
Similarly, consider $G/N$ space $(X/N)^{HN/N}$ and its invariant
subspace $(A/N)^{HN/N}$ for $H\in S$. It is easy to see that
$I_{f_{1}}^{HN/N}=J_{f_{1}}^{HN/N}$ where
\begin{align*}
J_{f_{1}}^{HN/N}&=\{a\in H^{*}(B_{G/N})\hbox{\rm :}\ a\otimes
f_{1}\in \hbox{Im}\,\{H_{G/N}^{*}((A/N)^{HN/N})\\[.4pc]
&\quad\,\rightarrow H_{G/N}^{*}((A/N)^{G/N})\}\}.
\end{align*}

This shows that if $X/N$ is totally nonhomologous to zero in
$(X/N)_{G/N}\rightarrow B_{G/N}$, then the local weights of $A$
are equal to the local weights of $X$.\hfill \epr
\end{proof}

\section*{Acknowledgements}

The authors would like to thank the referee for his suggestions
and corrections.

\end{document}